\documentclass[10pt]{amsart}
\usepackage{amsmath}
\usepackage{amssymb}
\usepackage{enumerate}
\usepackage{amsbsy}
\usepackage{amsfonts}
\usepackage{color}

\newtheorem{thm}{Theorem}[section]
\newtheorem{prop}[thm]{Proposition}
\newtheorem{lem}[thm]{Lemma}

\newtheorem{cor}[thm]{Corollary}
\newtheorem{defn}[thm]{Definition}

\numberwithin{equation}{section}

\theoremstyle{remark}
\newtheorem{rem}{Remark}

\newcommand{\supp}{\text{supp }}

\newcommand{\R}{\mathbb{R}}
\newcommand{\C}{\mathbb{C}}
\newcommand{\ep}{\varepsilon}

\newcommand{\al}{\alpha}

\newcommand{\ls}{\lesssim}
\newcommand{\ld}{\lambda}

\newcommand{\la}{\langle}
\newcommand{\ra}{\rangle}

\begin{document}

\newcommand{\sar}{S_\al(\mathbb R)}
\newcommand{\xar}{X_\al(\mathbb R)}
\newcommand{\yar}{Y_\al(\mathbb R)}

\newcommand{\ald}{\mathcal A}

\newcommand{\blb}{\big[}
\newcommand{\brb}{\big]}
\newcommand{\re}{\mbox{Re}}
\newcommand{\mlf}{\mathbf{m}_{1,\ld}}
\newcommand{\mlft}{\widetilde{\mathbf{m}}_{1, \ld}}
\newcommand{\mls}{\mathbf{m}_{2, \ld}}
\newcommand{\im}{\mbox{Im}}
\newcommand{\bl}{\big\langle}
\newcommand{\br}{\big\rangle}
\newcommand{\ar}{\mathcal A_\lambda}
\newcommand{\da}{{\delta(2-\alpha)}}
\newcommand{\sa}{{\mathcal S_\alpha}}

\newcommand{\uoq}{\mathbf u_{1, R_1}}
\newcommand{\utq}{\mathbf u_{2, R_2}}
\newcommand{\ui}{\mathbf u_{i}}
\newcommand{\uon}{\mathbf u_{1}}
\newcommand{\utw}{\mathbf u_{2}}
\newcommand{\uth}{\mathbf u_{3}}
\newcommand{\vn}{\mathbf v}

\newcommand{\ul}{U_{m}}
\newcommand{\ult}{U_{m}^2}
\newcommand{\vl}{V_{m}}

\title[Scattering for semirelativistic Hartree equations]
      {Small data scattering of semirelativistic Hartree equation}

\author[C. Yang]{Changhun Yang}
\address{Department of Mathematical Sciences, Seoul National University, Seoul 151-747, Republic of Korea}
\email{maticionych@snu.ac.kr}

\begin{abstract}
In this paper we study the small data scattering of Hartree type semirelativistic equation in space dimension $3$. The Hartree type nonlinearity is $[V * |u|^2]u$ and the potential $V$ which generalizes the Yukawa has some growth condition. We show that the solution scatters to linear solution if an initial data given in $ H^{s,1}$ is sufficiently small and $s>\frac14$. Here, $H^{s, 1}$ is  Sobolev type space taking in angular regularity with norm defined by $\|\varphi\|_{ H^{s, 1}} = \|\varphi\|_{ H^{s}} + \|\nabla_{\mathbb S} \varphi\|_{H^{s}}$. To establish the results  we employ the recently developed Strichartz estimate which is $L_\theta^2$-averaged on the unit sphere $\mathbb S^{2}$ and construct the resolution space based on $U^p$-$V^p$ space.
\end{abstract}

\thanks{2010 {\it Mathematics Subject Classification.} 35Q55, 35Q53. }
\thanks{{\it Key words and phrases.} semirelativistic Hartree equation, Yukawa type potential, small data scattering, angularly averaged Strichartz estimate, $U^p$ and $V^p$ spaces}

\maketitle

\newcommand{\na}{|\nabla|^\al}
\section{Introduction}
In this paper we consider the following Cauchy problem:
\begin{align}\label{main eqn}\left\{\begin{array}{l}
i\partial_tu = \Lambda_m u + F(u)\;\;
\mbox{in}\;\;\mathbb{R}^{1+3},
\\
u(x,0) = \varphi(x) \;\; \mbox{in}\;\; \mathbb{R}^3,
\end{array}\right.
\end{align}
where $\Lambda_m$ is the fourier multiplier defined by $\Lambda_m=(m-\Delta)^\frac 12$ and $F(u)$ is nonlinear term of Hartree type such that $F(u) = [V * |u|^2]u$ with a smooth $V$ in $\mathbb R^3 \setminus \{0\}$. Here $m>0$ is mass and $*$ denotes the convolution in $\R^3$. The concerned Hartree potential is defined as follows:
\begin{defn}
For $0\le \gamma_1,\gamma_2<3$ the potential $V$ is said to be of type $(\gamma_1,\gamma_2)$ if it satisfies the growth condition such that $\widehat{V} \in C^{4}(\mathbb R^3 \setminus\{0\})$ and  for $ 0 \le k \le 4$ \begin{align}\label{growth}
|{\nabla^k}\widehat V(\xi)| \lesssim |\xi|^{-\gamma_1-k}\;\;\mbox{for}\;\; |\xi| \le 1, 
\quad |{\nabla^k}\widehat V(\xi)| \lesssim |\xi|^{-\gamma_2-k}\;\;\mbox{for}\;\;|\xi| > 1.
\end{align}
\end{defn}
The Coulomb potential $V(x)=|x|^{-1}$ is of such type corresponding to $\gamma_1=\gamma_2=2$ and
the Yukawa potential $V(x)=e^{-\mu_0 |x|} |x|^{-1}, \ \mu_0>0$ is  corresponding to $\gamma_1=0, \gamma_2=2$.
The equation \eqref{main eqn} with these two potentials, which is called semirelativistic Hartree equation, arises in the mean-field limit of large systems of bosons, see, e.g., \cite{E07,l03,EH87}. 
In this paper we study \eqref{main eqn} with the above generalized potentials.

By Duhamel's formula, \eqref{main eqn} is written as an
integral equation
\begin{equation}\label{integral}
u = e^{-it{\Lambda_m}}\varphi -i \int_0^t e^{-i(t-s){\Lambda_m}}F(u)(s)\,ds.
\end{equation}
Here we define the linear propagator $e^{-it{\Lambda_m}}$ given by the solution to the linear problem $i\partial_t v =  {\Lambda_m} v$ with initial datum $v(0) = \varphi$. It is
formally written  by
\begin{align}\label{int eqn}
e^{-it{\Lambda_m}}\varphi = \mathcal F^{-1} \big(e^{-it\sqrt{m +|\xi|^2}} \mathcal F(\varphi)\big) = (2\pi)^{-3}\int_{\mathbb{R}^3} e^{i( x\cdot \xi - t\sqrt{m + |\xi|^2} )}\widehat{\varphi}(\xi)\,d\xi.
\end{align}

The purpose of this research is to study the existence and uniqueness of solutions and observe the behaviour of solutions as time goes to infinity, in particular  comparing them with the linear solutions, which is well-known in area of dispersive PDE as the well-posedness and scattering problem respectively. The following is formal definition of scattering. 
\begin{defn}\label{scattering defn}
	We say that a solution $u$ to \eqref{main eqn} scatters (to $u_\pm$) in a Hilbert space $\mathcal H$ if there exist  $\varphi_\pm \in \mathcal H$ (with $u_\pm(t) = e^{-it{\Lambda_m}}\varphi_\pm$) such that $\lim_{t  \to \pm\infty}\|u(t) - u_\pm\|_{\mathcal H} = 0$.
\end{defn}
One of candidates for Hilbert space $\mathcal H$ is the Sobolev spaces $H^s(\R^3)$. That is, we will show the well-posedness and scattering results when the initial data is given in $H^s(\R^3)$. Especially we want to find the minimum value of $s$ that ensures the scattering states of corresponding solutions, which is called low regularity problem.

There have been a lot of results on this subject.
Firstly, Lenzmann \cite{l07} established the global existence of solutions for Yukawa type potential using energy methods provided the initial data given in $H^{\frac12}(\R^3)$ is sufficiently small.
Herr and Lenzmann \cite{HL2014} showed that for Coulomb type potential the almost optimal local well-posedness holds for initial data with $s>\frac14$ (and $s>0$ if the data is radially symmetric) using localized Strichartz estimates and the Bourgain spaces.
In \cite{chonak,choz0,co07,coss 09}, they considered the generalized potential from Coulomb type, namely, $V(x)=|x|^{-\gamma}$, for $0<\gamma<3$ (corresponding to $\gamma_1=\gamma_2=3-\gamma$ in our definition) and investigated well-posedness and scattering of equations.
The most recent results on the Yukawa potential were obtained by Herr and Tesfahun \cite{hete} where they showed the small data scattering result for initial data with $s>\frac12$( and $s>0$ if the data is radially symmetric) using $U^p-V^p$ spaces method which has proved effective to derive scattering result.

In this paper we consider the range $0<s\le \frac12$ where the scattering result has been proved only when radial assumption is given to initial data \cite{hete} and aim to obtain the similar result with a weaker assumption.  
We prove the scattering result when $s>\frac14$ by imposing additional one angular regularity to the initial data. 
Let us introduce angular derivative and angularly regular Sobolev space. The spherical gradient $\nabla_{\mathbb S}$ is restriction of the gradient on the unit sphere which is well-defined, that is, independent of coordinates of $\mathbb{S}^2$. It satisfies a following relation
$$\nabla=\frac1r\nabla_{\mathbb S}+\theta\frac{\partial}{\partial r}, \ x=r\theta, \ \theta\in\mathbb{S}^2,$$
and also has a concrete formula $\nabla_{\mathbb S}=x \times \nabla$. 
A function space $H^{s, 1}$ is the set of all $ H^s$ functions whose angular derivative is also in $ H^s$. The norm is defined by $\|f\|_{ H^{s, 1}} := \|f\|_{ H^s} + \|\nabla_{\mathbb S} f\|_{ H^{s}}$. It contains all radially symmetric functions. 

Our main result is the following.
\begin{thm}\label{main thm}
Let $s>\frac14$. Suppose the potential $V$ in \eqref{main eqn} is radially symmetric and of type $(\gamma_1,\gamma_2)$ with $0\le\gamma_1<1$ and $\frac32<\gamma_2<3$. Then there exists $\delta > 0$ such that for any $\varphi \in H^{s, 1}$ with $\|\varphi\|_{ H^{s, 1}} \le \delta$, \eqref{main eqn} has a unique solution $u \in (C \cap L^\infty)(\mathbb R; H^{s, 1})$ which scatters in $ H^{s, 1}$.
\end{thm}

\begin{rem}
The potential in Theorem~\ref{main thm} includes the Yukawa. Concerning the Coulomb potential, non-existence scattering results \cite{choz0} and modified scattering results \cite{F2014} have been established.
\end{rem}

Our proof is fundamentally based on fixed point argument and Littlewood-Paley decomposition. In order to occur a contraction, we use  frequency-localized spherical Strichartz estimates and construct a resolution space using $U^p, V^p$ spaces where linear estimates for free solutions could be transferred.



We have an application to the following Hartree Dirac equations:
\begin{align}\label{sub eqn}\left\{\begin{array}{l}
i\partial_t\psi = ({\alpha}\cdot D+m\beta) \psi + [V * |\psi|^2]\psi \;\;
\mbox{in}\;\;\mathbb{R}^{1+3},
\\
\psi(x,0) = \psi_0(x) \;\;x \in \mathbb{R}^3,
\end{array}\right.
\end{align}
where $D=-i\nabla, \psi :\R^{1+3}\rightarrow\mathbb{C}^4$ is the Dirac spinor, $m>0$ is mass and $\beta$ and $\alpha=(\alpha_1,\alpha_2,\alpha_3)$ are the Dirac matrices.  
If $V$ is Coulomb potential, \eqref{sub eqn} appears when Maxwell-Dirac system with zero magnetic field is uncoupled \cite{cg 76}. 
And in the same paper \cite{cg 76} it is conjectured that \eqref{sub eqn} with Yukawa also might be obtained by uncoupling Dirac-Klein-Gordon system as Maxwell-Dirac case.
For more information about Dirac equation, see \cite{cg 76} and references therein. 

Following \cite{dfs07} (see also \cite{BH2017}) we introduce the projection operators $\Pi_{\pm}^m(D)$ with symbol 
\begin{equation*}
\Pi_{\pm}^m(\xi)=\frac12[I\pm\frac{1}{\langle\xi\rangle}_m(\xi\cdot\alpha+m\beta)].
\end{equation*} 
We then define $\psi_{\pm}:=\Pi_{\pm}^m(D)\psi$ and split $\psi=\psi_++\psi_-$.
By applying the operators $\Pi_{\pm}^m(D)$ to the equation \eqref{sub eqn}, and using the identity 
$$ \alpha\cdot D+m\beta=\Lambda_m(\Pi_{+}^m(D)-\Pi_{-}^m(D)) $$
we obtain
the following system of equations
\begin{equation}\label{system}
\begin{cases}
&(-i\partial_t+\Lambda_m)\psi_{+}= \Pi_{+}^m(D)[(V*|\psi|^2)\psi], \\
&(-i\partial_t-\Lambda_m)\psi_{-}=\Pi_{-}^m(D)[(V*|\psi|^2)\psi], 
\end{cases}
\end{equation}
with initial data $\psi_0^{\pm}=\Pi_{\pm}^m(D)\psi_0$.
Observe that the linear propagators for this system have same formula as \eqref{main eqn} except for sign and the nonlinear term is also same except for projection operator. Note that Strichartz estimates holds regardless of sign and the Sobolev norm has an equivalence under this operator, i.e., 
$\| \Pi_{+}^m(D) f\|_{H^s} + \| \Pi_{+}^m(D) f\|_{H^s} \sim \| f\|_{H^s}.$
Since our proof for Theorem~\ref{main thm} does not require any structure of equation but relies on Strichartz estimates, function spaces and Littelwood-Paley decomposition, one can easily check the following Corollary:
\begin{cor}
	Let $s>\frac14$. Suppose the potential $V$ in \eqref{main eqn} is radially symmetric and of type $(\gamma_1,\gamma_2)$ with $0\le\gamma_1<1$ and $\frac32<\gamma_2<3$. Then there exists $\delta > 0$ such that for any $\psi_0 \in H^{s, 1}$ with $\|\psi_0\|_{\dot H^{s, 1}} \le \delta$, \eqref{sub eqn} has a unique solution $\psi \in (C \cap L^\infty)(\mathbb R; H^{s, 1})$ which scatters in $ H^{s, 1}$.
\end{cor}

\medskip

\section{Notations and Preliminaries}
\subsection{Notations}
The Fourier transform of $f$ is denoted by $\widehat{f} = \mathcal F( f)$  and the inverse Fourier transform is by $\mathcal F^{-1}$ such that
$$\mathcal{F}(f)(\xi) = \int_{\mathbb{R}^3} e^{- ix\cdot \xi} f(x)\,dx,\quad \mathcal F^{-1} (g)(x) = (2\pi)^{-3}\int_{\mathbb{R}^3} e^{ix\cdot \xi} g(\xi)\,d\xi.$$ 
We denote the frequency variables by capital letters $M,N>0$ which is assumed dyadic number, that is of the form $2^m,2^n$ with $m,n\in\mathbb{Z}$.
Let $\rho\in C_{0,rad}^\infty(-2,2)$ be such that $\rho(s)=1$ if $|s|<1$ and $\chi_M$ be defined by $\chi_M(s)=\rho(\frac{s}{M})-\rho(\frac{2s}{M})$ for $M>0$. Then $\supp\chi_{M}=\{ s\in\R: \frac M2<|s|<2M \}$.
Fix $N_0\gg1$ and let $\beta_{N_0}:=\sum_{M \le N_0}\chi_{M}$ and $\beta_N:=\chi_N$ for $N>N_0$. Then $\supp\beta_{N_0}=\{ s\in\R: |s|<2N_0 \}$ and $\supp\beta_{N}=\supp\chi_N$ for $N>N_0$.
Denote $\widetilde\chi_M(s):= \chi_M(s/2) + \chi_M(s) + \chi_M(2s)$ and $\widetilde \beta_N$ similarly.
Next we define the Littlewood-Paley operators $\dot P_M$ and $P_N$ by $\mathcal F (\dot P_Mf) = \chi_M\widehat{f}$ for $M>0$ and $\mathcal F ( P_Nf) = \beta_N\widehat{f}$ for $N\ge N_0$ respectively. Further, define similarly $\widetilde {\dot P_M} $ and $\widetilde P_N$ using $\widetilde\chi_M$ and $\widetilde\beta_N$. Then $\widetilde{ \dot P_M }\dot {P}_M = \dot P_M \widetilde{ \dot P_M }= \dot P_M$ and $\widetilde P_N P_N=P_N \widetilde P_N = P_N$.

We denote $L^r = L_x^r(\mathbb R^3)$ and  $\mathcal{L}^r=\mathcal{L}_\rho^r(\R)=L_\rho^r(\rho^{2}d\rho)$ for  $1 \le r \le \infty$. Consider the mixed-normed space. For a Banach space $X$, $u \in L_I^q X$ iff $u(t) \in X$ for a.e. $t \in I$ and $\|u\|_{L_I^qX} := \|\|u(t)\|_X\|_{L_I^q} < \infty$. We denote  $L_I^qX = L_t^q(I; X)$ and   $L_t^qX = L_{\mathbb R}^qX$. 

Positive constants depending only on $m,N_0$ are denoted by the same letter $C$, if not specified. $A \lesssim B$ and $A \gtrsim B$ means that $A \le CB$ and
$A \ge C^{-1}B$, respectively for some $C>0$. $A \sim B$ means that $A \lesssim B$ and $A \gtrsim B$.

\subsection{Function spaces}
In this subsection we introduce the $U^p, V^p$ function spaces. For the general theory, see e.g. \cite{hhk}, \cite{ktv}.

Let $1\leq p<\infty$. We call a finite set $\{t_0,\ldots,t_J\}$ a partition if $-\infty<t_0<t_1<\ldots<t_J\leq \infty$, and denote the set of all partitions by $\mathcal{T}$. A corresponding step-function $a:\R\to L^2(\R^3)$ is called $U^p$-atom if
\[
a(t)=\sum_{j=1}^J \mathbf{1}_{[t_{j-1},t_j)} (t) f_j, \quad \sum_{j=1}^J\|f_j\|_{L^2(\R^3)}^p=1,\quad \{t_0,\ldots,t_J\}\in \mathcal{T},
\]
and $U^p$ is the atomic space. The norm is defined by
\begin{align*}
\| u\|_{U^p} :=
\inf \Big\{ \sum_{k=1}^\infty |\lambda_k| : u=\sum_{k=1}^\infty\lambda_k a_k, \text{ where } a_k \text{ are } U^p\text{-atoms and } \lambda_k\in\C \Big\}.
\end{align*}
Further, let $V^p$ be the space of all right-continuous $v:\R \to L^2(\R^3)$ satisfying
\begin{align}\label{V norm}
\|v\|_{V^p}:=\sup_{\{t_0,\ldots,t_J\}\in \mathcal{T}}\big(\sum_{j=1}^J\|v(t_j)-v(t_{j-1})\|_{L^2(\R^3)}^p\big)^{\frac1p}.
\end{align}
with the convention $v(t_J)=0$ if $t_J=\infty$. 
Likewise, let $V_-^p$ denote the spaces of all functions $v:\R\rightarrow L^2(\R^3)$ satisfying $v(-\infty)=0$ and $\|v\|_{V^p}<\infty$, equipped with the norm \eqref{V norm}.
We define $V_{-,rc}^p$ by the closed subspace of all right continuous $V_{-}^p$ functions. 

Now we list some useful Lemmas on $U^p,V^p$ spaces.

\begin{lem}\label{Lem:property UV}
Let $1\le p<q<\infty$.
\begin{enumerate}
\item $U^p,V^p,V_-^p$ and $V_{-,rc}^p$ is Banach spaces.
\item The embeddings $U^p\hookrightarrow V^p_{-,rc} \hookrightarrow U^q \hookrightarrow L^\infty(\R;L^2)$ are continuous.
\item The embeddings $V^p\hookrightarrow V^q$ and $V_-^p\hookrightarrow V_-^q$ are continuous.
\item (Duality) For $1<p<\infty$, $\| u\|_{U^p}= \sup_{\{v\in V^{p'}:\|v\|_{V^{p'}=1}\}} \int_{-\infty}^{\infty}\la u'(t), v(t)\ra_{L_x^2} dt.$
\end{enumerate}
\end{lem}

\begin{defn}[Adapted function spaces]
We define $U_m^p$(and $V_m^p$ respectively) by the spaces of all functions $u$ such that $e^{it\Lambda_m}u\in U^p$ ($e^{it\Lambda_m}v\in V^p$ respectively) with the norm
\begin{align*}
\| u\|_{U_m^p}:=\|e^{it\Lambda_m} u \|_{U^p} \quad
( \ \| v\|_{V_m^p}:= \|e^{it\Lambda_m} v \|_{V^p} \text{ respectively }).
\end{align*}
\end{defn}
The properties in Lemma~\ref{Lem:property UV} also hold for the spaces $U_m^p$ and $V_m^p$.

\begin{lem}[Transfer principle]\label{Lem:transfer}
	Let $T: L^2 \to L_{loc}^1(\mathbb R^3; \mathbb C)$ be a linear operator satisfying that
	$$
	\|T(e^{-it\Lambda_m} f )\|_{L_t^q X} \lesssim \| f \|_{L^2}
	$$
	for some $1 \le q < \infty$ and a Banach space $X \subset L_{loc}^1(\R^3;\C)$. Then
	$$
	\|T(u)\|_{L_t^qX} \lesssim \|u\|_{{U_m^q}}.
	$$
\end{lem}

\subsection{Strichartz estimates}
Let the pair $(q, r)$ satisfy that $2 \le q, r \le \infty$, $\frac 2{q} + \frac 3{r} = \frac 32$.
Then it holds from \cite{cox}
\begin{align}
\bl M \br^{-\frac{5}{3q}}\|e^{-it{\Lambda}_m} \dot P_M\varphi\|_{L_t^{q} L_x^r} \lesssim \| \dot P_M\varphi\|_{L_x^2}.\label{besov str}
\end{align}
The case $q=2$ is sufficient in our discussion. For the $N_0 \gg 1$ we have 
\begin{align*}
\|e^{-it\Lambda_m}P_{N_0}\varphi\|_{L_t^2L_x^6} \lesssim_{N_0} \|P_{N_0}\varphi\|_{L^2},
\end{align*}
which gives by transfer principle in Lemma~\ref{Lem:transfer}
\begin{align}\label{low str}
\|P_{N_0}u\|_{L_t^2L_x^6} \lesssim_{N_0} \|P_{N_0}\varphi\|_{U_m^2}.
\end{align}
This endpoint estimate can be extended to a wider range with weaker angular integrability in the left term. That is, we consider the following $L_t^q\mathcal L_\rho^r L_\theta^{{r_*}}$ norm
with $r_*\le r < \infty$  defined by
$$
\|u\|_{L_t^q \mathcal L_\rho^r L_\theta^{{r_*}}} = \left(\int_\mathbb R \left\|\big(\int_{\mathbb S^2}|u(t, \rho \theta)|^{{r_*}}\,d\theta\big)^\frac1{{r_*}}\right\|_{L_\rho^r(\rho^2 d\rho)}^q dt \right)^\frac1q.
$$
If $r = \infty$, then we define $\mathcal L_\rho^\infty = L_\rho^\infty$.
Then for $\frac{10}{3} < r < 6$, there holds 
\begin{align*}
\|e^{-it{\Lambda_m}}  \dot P_M\varphi\|_{L_t^{2} \mathcal L_\rho^r L_\theta^2}  \lesssim \| \dot P_M\varphi\|_{L_x^2} \times
\begin{cases}
M^{\frac12 - \frac 3r}\bl M\br^{-\frac12+\frac4r} ,\  &\text{if} \ \frac{10}{3} < r <  4\\
M^{-\frac14}\bl M\br^{\frac12+\ep}, \ &\text{if} \ r=4 \ \text{and}\ \epsilon>0, \\
M^{1-\frac3r}, \ &\text{if} \ 4<r<6. \end{cases}
\end{align*}
For this see the Klein-Gordon case of Theorem 3.3 in \cite{ghn}.
Especially if $\frac{10}{3} < r <  4$ and $N>N_0$ we have for $u\in U_m^2$ by transfer principle into $U_m^2$ spaces
\begin{align}\label{ang str}
\|  P_N u\|_{L_t^{2} \mathcal L_\rho^r L_\theta^2}  \lesssim N^\frac1r \| P_N u\|_{U_m^2},
\end{align}
which we will intensively use in following argument.
\begin{rem}
	Note that the minimum loss of regularity occurs when $r$ is close to $4$. And this is essentially related with the regularity condition on initial data $s>\frac14$. It can be easily checked that in the range $4<r<6$ the bound is sharp if we consider the homogeneous case and scaling argument, but in the other range the sharpness is not known yet. If we can improve the bound in this range we might obtain better regularity result, i.e., threshold of well-posedness could be lowered.  	
\end{rem}

\subsection{Properties of angular derivative}
In this section we introduce a series of lemmas concerning angular derivative.
\begin{lem}\label{radial1}
Let $\psi, f$ be smooth and let $\psi$  be radially symmetric. Then $$\nabla_{\mathbb{S}} (\psi * f) = \psi * \nabla_{\mathbb{S}} f.$$
\end{lem}
From this we check the order of the projection operator and angular derivative can be reversed: $\nabla_{\mathbb S} \dot P_M f=\dot P_M \nabla_{\mathbb S} f$ for $M>0$.

The next one is on the Sobolev inequality on the unit sphere \cite{crt-n}.
\begin{lem}\label{sob-s2}
For any $2 < \widetilde r < \infty$
$$
\|f\|_{L_\theta^{\widetilde r}(\mathbb{S}^2)} \lesssim \|f\|_{L_\theta^2(\mathbb{S}^2)} + \|\nabla_{\mathbb S}f\|_{L_\theta^2(\mathbb{S}^2)}, \quad 
\|f\|_{L_\theta^\infty(\mathbb{S}^2)} \lesssim \|f\|_{L_\theta^{\widetilde r}(\mathbb{S}^2)} + \|\nabla_{\mathbb S}f\|_{L_\theta^{\widetilde r}(\mathbb{S}^2)}.
$$
\end{lem}

The final one is extended Young's convolution estimates.
\begin{lem}[Lemma 7.1 of \cite{chonak}]\label{radial2}
	If $\psi$ is radially symmetric, then
	$$
	\|\psi*f\|_{\mathcal L_\rho^p L^q_\theta} 
	\le \|\psi\|_{L^{p_2}_x} \|f\|_{\mathcal L^{p_1}_\rho L^{q_1}_\theta} ,
	$$
	for all $p_1,p_2,p,q,q_1 \in [1,\infty]$ satisfying
	$$
	\frac 1{p_1} + \frac 1{p_2} - 1 = \frac 1p, \quad \frac 1{q_1} + \frac 1{p_2} - 1 \le \frac 1q.$$
\end{lem}

\subsection{Norm of Potential}
We calculate the $L^p$ norm of $\dot P_M V$ and $\mathcal{F}^{-1} \chi_M$ for $1<p<\infty$. We simply denote $\dot P_M V$ by $V_M$.
\begin{align*}
\int |V_M(x)|^pdx &= \int_{|x| \le M^{-1}} |V_M(x)|^pdx + \int_{|x| > M^{-1}}|x|^{-4p}|x|^{4p} |V_M(x)|^p dx  \\
&\lesssim M^{-3}\|V_M\|_{L^\infty}^p + M^{4p-3}\||x|^{4}V_M\|_{L^\infty}^p\\
&\lesssim M^{-3} \| \chi_M(\xi)\widehat V(\xi) \|_{L^1}^p + M^{4p-3}\|\nabla_\xi^{4}\big( \chi_M(\xi) \widehat V(\xi) \big)\|_{L^1}^p.
\end{align*}
Using the assumption \eqref{growth} of $V$ we estimate $ \|\nabla_\xi^{k}\big( \chi_M(\xi) \widehat V(\xi) \big)\|_{L^1} \ls M^{-k-\gamma+3}$ for $0\le k\le 4$, where $\gamma=\gamma_1$ if $0<M\le1$, or 
$\gamma=\gamma_2$ if $ M>1$.
Thus we have
\begin{equation}\label{ineq:potential}
\|V_M \|_{L_x^p} \lesssim 
\begin{cases}
M^{3-\frac3p-\gamma_1}\   ,&\text{if} \  0<M\le1 \\ 
M^{3-\frac3p-\gamma_2}\   ,&\text{if} \ M>1. \end{cases}
\end{equation}
Also we can check by simple calculation 
\begin{align*}
\| \mathcal{F}^{-1} \chi_M \|_{L_x^{p}} \ls M^{3-\frac3p}.
\end{align*} 
 
Now we are ready to prove the main theorem.

\section{Proof of Main theorem}\label{main}

Let us define the  Banach space $ X^s$ by
$$
X^s := \Big\{u : \mathbb R \to  H^{s} \Big|\; P_Nu, \nabla_{\mathbb{S}} P_Nu \in \ul^2(\mathbb R; L_x^2) \;\;\forall N \ge N_0  \Big\}
$$
with the norm  
$$\|u\|_{ X^s} = \left(  \sum_{N \ge N_0}  { N }^{2s}\| P_N u\|_{\ul^{2,1}}^2\right)^\frac12, 
\;\mbox{where}\; \|u\|_{\ul^{2,1}} = \|u\|_{\ul^2} + \|\nabla_{\mathbb S}u\|_{\ul^2}.$$
Let $X_+^s$ be the restricted space defined by
$$
 X_+^s = \Big\{ u \in C([0, \infty);  H^{s}) \Big|\; \chi_{[0, \infty)}(t)u(t) \in X^s\Big\}
$$
with norm $\|u\|_{X_+^s} := \|\chi_{[0, \infty)}u\|_{X^s}$.

Let $\mathcal D_+^s(\delta)$ be a complete metric space $\{u \in X_+^s \big|\; \|u\|_{X_+^s} \le \delta\}$ equipped with the metric $d(u, v) := \|u-v\|_{X_+^s}$.
Then we will show that the nonlinear functional $\Psi(u) = e^{-it\Lambda_m}\varphi + \mathcal N_m(u,u,u)$ is a contraction on $D_+^s(\delta)$, where
$$
\mathcal N_m (u_1,u_2,u_3) = -i\int_0^t e^{-i(t-t')\Lambda_m} [V*(u_1\bar{u}_2)u_3]\,dt'.
$$
Clearly, $\|e^{-it\Lambda_m} \varphi\|_{X_+^s} \lesssim \|\varphi\|_{H^{s, 1}}$ so  it suffices to show that
\begin{align}\label{contract}
\|\mathcal N_m(u_1,u_2,u_3)\|_{X_+^s} \lesssim \prod_{j=1}^3\|u_j\|_{X_+^s}^3
\end{align}
This readily implies estimates for difference 
$$ \|\mathcal N_m(u,u,u) - \mathcal N_m(v,v,v)\|_{X_+^s} \lesssim (\|u\|_{X_+^s} + \|v\|_{X_+^s})^2\|u-v\|_{X_+^s}, \ \text{for} \ u,v\in X^s$$
and thus we can find $\delta$ small enough for $\Psi$ to be a contraction mapping on $\mathcal D_+^s(\delta)$. From now, we simply denote $N_m(u,u,u)$ by $N_m(u)$.

Since $e^{it\Lambda_m}P_N \mathcal N_m(u)$ and $e^{it\Lambda_m}P_N \nabla_{\mathbb S}\mathcal N_m(u)$ are in $V_{-, rc}^2(\mathbb R; L^2)$, and
$$
\sum_{N \ge N_0} { N }^{2s}(\|e^{it\Lambda_m}P_N \mathcal N_m(u)\|_{V^2} + \|e^{it\Lambda_m}P_N \nabla_{\mathbb S}\mathcal N_m(u)\|_{V^2})^2 < \infty
$$
from \eqref{contract}, $\lim_{t \to +\infty} e^{it\Lambda_m}\mathcal N_m(u)$ exists in $H^{s, 1}$. 
Define a scattering state $u_+$ with $$\varphi_+:= \varphi + \lim_{t\to +\infty} e^{it\Lambda_m}\mathcal N_m(u).$$ By time symmetry we can argue in a similar way for the negative time. Thus we get the desired result.

We start to show \eqref{contract}. We may assume that $u(t) = 0$ for $-\infty < t < 0$.
From the duality in Lemma~\ref{Lem:property UV},
\begin{align*}
\|P_N \mathcal N_m(u_1,u_2,u_3)\|_{\ul^2}
\ls \sup_{\|v\|_{\vl^2} \le 1}\Big|\int_{\mathbb R}\int_{\mathbb R^3} [V * (u_1\bar{u}_2)]u_3(t)P_N\overline{v(t)}\,dxdt\Big|.
\end{align*}
Using Littlewood-Paley decomposition and applying Lemma~\ref{radial1} and Leibniz rule, we have
\begin{align}\label{cubicbound0}
\|\mathcal N_m(u)\|_{X_+^s}^2
\le \sum_{N \ge N_0}{ N }^{2s} \sup_{\|v\|_{\vl^2} \le 1}\left(\sum_{N_1, N_2, N_3 \ge N_0}
\sum_{k=0}^3 I_k(N,N_1,N_2,N_3)
\right)^2, 
\end{align}
where 
\begin{align*}
I_0(N_1,N_2,N_3,N)&=\bigg|\iint V *(P_{N_1}u_1 P_{N_2}\bar{u}_2)P_{N_3}u_3 P_N \bar{v}\,dxdt \bigg|, \\
I_1(N_1,N_2,N_3,N)&=\bigg|\iint V *(\nabla_{\mathbb S} P_{N_1}u_1 P_{N_2}\bar{u}_2)P_{N_3}u_3 P_N \bar{v}\,dxdt \bigg|, \\
I_2(N_1,N_2,N_3,N)&=\bigg|\iint V *(P_{N_1}u_1 \nabla_{\mathbb S} P_{N_2}\bar{u}_2)P_{N_3}u_3 P_N \bar{v}\,dxdt \bigg|, \\
I_3(N_1,N_2,N_3,N)&=\bigg|\iint V *(P_{N_1}u_1 P_{N_2}\bar{u}_2)\nabla_{\mathbb S}
 P_{N_3}u_3 P_N \bar{v}\,dxdt \bigg|.
\end{align*}
Since the argument will not be affected by complex conjugation, we drop the conjugate symbol. By Lemma~\ref{radial1} we can change the order of deriavtive operator $\nabla_{\mathbb S}$ and projection $P$. 
Thus to show \eqref{contract}, we suffices to prove the following:
\begin{align}\label{ineq:goal}
\sum_{N \ge N_0}{ N }^{2s} \sup_{\|v\|_{\vl^2} \le 1}\left(\sum_{N_1, N_2, N_3 \ge N_0}
I(N,N_1,N_2,N_3)
\right)^2
\ls
\prod_{i = 1, 2, 3}\|u_i\|_{X_+^s}^2,
\end{align}
where 
$$I(N_1,N_2,N_3,N) := \Big|\iint V *(P_{N_1}\uon P_{N_2}\utw)P_{N_3}\uth P_N v\,dxdt \Big|,$$
and at most one of $\mathbf{u_i}$ could take an angular derivative $\nabla_{\mathbb S}$, i.e. $\mathbf{u_i}= \nabla_{\mathbb S} u_i$. To prove this inequality we introduce the following proposition.
\begin{prop}\label{prop}
Let $s>\frac14$. Suppose $P_{N_i}u_i\in\ul^{2,1}$, $P_{N}v \in\vl^{2}$ for $i=1,2,3$
and at most one of $\mathbf{u_i}$ could take an angular derivative $\nabla_{\mathbb S}$. Then for $\frac14<\frac1r<\min(s,\frac{\gamma_2}{6},\frac{3}{10})$ it holds 
\begin{equation}\label{ineq:key}
\begin{split}
 I(N_1,N_2,N_3,N) &\ls C(N,N_1,N_2,N_3)
\|P_{N_1}u_1\|_{\ul^{2,1} }  \|P_{N_2}u_2\|_{\ul^{2,1} } 
\|P_{N_3}u_3  \|_{\ul^{2,1}} \|P_N v \|_{\vl^2} , \\
C(N,N_1,N_2,N_3) &= \begin{cases}
N_1^\frac1r N_2^\frac1r& \ \text{for} \ N_3\gtrsim N, \\
\min(N_1,N_2)^\frac1r N_3^\frac1r & \ N_3\ll N.  
\end{cases}  
\end{split}
\end{equation}
Here  the implicit constant only depends on $r, N_0$.
\end{prop}

Now, we postpone the proof of Proposition~\ref{prop} in a moment and explain how this result implies \eqref{ineq:goal}.
Let us split the summation of LHS in \eqref{ineq:goal}  into two parts as follows:
$$
\mbox{LHS of}\;\; \eqref{ineq:goal} = \sum_{N_3\gtrsim N}+\sum_{N_3\ll N}:=S_1+S_2.
$$
Fix $r$ as in Proposition~\ref{prop}. Apply the first case of \eqref{ineq:key} to $S_1$
\begin{align*}
S_1
&\lesssim
\sum_{N \ge N_0 } { N }^{2s}
\left(
\sum_{N_1, N_2\ge N_0} N_1^\frac1r \|P_{N_1}u_1\|_{\ul^{2,1} } N_2^\frac1r  \|P_{N_2}u_2\|_{\ul^{2,1} } 
\sum_{N_3 \gtrsim N}\|P_{N_3} u_3  \|_{\ul^{2,1}}  \right)^2\\
&\lesssim
\sum_{N \ge N_0 } 
\left(
\sum_{N_1, N_2\ge N_0} N_1^{\frac1r-s}N_2^{\frac1r-s}
N_1^{s}\|P_{N_1}u_1\|_{\ul^{2,1} } N_2^{s} \|P_{N_2}u_2\|_{\ul^{2,1} } 
\sum_{N_3 \gtrsim N} (\frac{N}{N_3})^s  N_3^s\|P_{N_3}u_3  \|_{\ul^{2,1}}  \right)^2\\
&\ls \|u_1\|_{X_+^s}^2\|u_2\|_{X_+^s}^2\sum_{N \ge N_0} \left(
\sum_{N_3\gtrsim N} (\frac{N}{N_3})^s
N_3^s \|P_{N_3}u_3\|_{\ul^{2,1} }  \right)^2 \\
&\lesssim
\prod_{i = 1, 2, 3}\|u_i\|_{X_+^s}^2.
\end{align*}
$S_2$ is estimated using  the second case of \eqref{ineq:key}. By symmetry we may assume $N_1\le N_2$.
\begin{align*}
S_2&\lesssim
\sum_{N \ge N_0} N^{2s} \left(
\sum_{\substack{N_0\le N_1 \le N_2\\ N_3\ll N}}     
N_1^{\frac1r-s}N_3^{\frac1r-s}
N_1^s\|P_{N_1}u_1\|_{\ul^{2,1} } \|P_{N_2}u_2\|_{\ul^{2,1} } 
N_3^s\|P_{N_3}u_3 \|_{\ul^{2,1}} \right)^2\\
&\ls \|u_1\|_{X_+^s}^2\|u_3\|_{X_+^s}^2\sum_{N \ge N_0} \left(
\sum_{N_2\gtrsim N} (\frac{N}{N_2})^s
N_2^s \|P_{N_2}u_2\|_{\ul^{2,1} }  \right)^2 \\
&\lesssim
\prod_{i = 1, 2, 3}\|u_i\|_{X_+^s}^2.
\end{align*}

So it remains to prove the Proposition~\ref{prop}. 
To simplify the notations, we assume all the functions are localized one, i.e., $P_{N_i}u_i=u_i$ for $i=1,2,3$ and $P_{N}v=v$. And we use the bold notation $\mathbf{u_i}$ when it could take an angular derivative or not. But be cautious that at most one of bold $\mathbf{u_i}$ could take. In other words, the estimates hold true even if at most one of $\mathbf{u_i}$ take an angular derivative.

\section{Proof of Proposition~\ref{prop}}
We perform an additional decomposition for $I$:
\begin{align*}
I(N_1,N_2,N_3,N)
&\lesssim  \sum_{M>0} |\iint \dot P_M V *(\uon\utw) \widetilde{\dot{P}_M}( \uth v) \,dxdt | \\
&=  \sum_{M>0} |\iint V_M *(\uon\utw) \widetilde{\dot{P}_M} ( \uth v) \,dxdt |,
\end{align*}
where at most one of bold $\mathbf{u_i}$ could take the angular derivative.
\subsection{Case1: $N_3\gtrsim N$}
In this subsection we prove that 
\begin{align*}
I\ls N_1^\frac1rN_2^\frac1r \|u_1\|_{\ul^{2,1} }  \|u_3\|_{\ul^{2,1} } 
\|u_3  \|_{\ul^{2,1}} \|v \|_{\vl^2}.
\end{align*}
Since the localized Strichartz estimates we apply have different admissible pair whether the support in frequency side is low part or not, that is, \eqref{low str} or \eqref{ang str}, we proceed to prove dividing the case whether $N_i$ is equal to $N_0$ or not for $i=1,2,3$.

Note that the support properties from Littlewood-Paley decomposition would restrict the range of summation over $M$.

\subsubsection{$N_0=\min(N_1,N_2)\sim \max(N_1,N_2)$}
In this case the support condition gives $M\ls N_0$. We estimate using H\"older and Young's inequality 
\begin{align*}
	I &\lesssim
	\sum_{M\ls N_0} \|V_M *(\uon \utw )\|_{L_t^1 L_x^\infty} 
	\| \widetilde{\dot{P}_M} (\uth v) \|_{L_t^\infty L_x^1} \\
	&\lesssim
	\sum_{M\ls N_0}   \| V_M \|_{L_x^\frac32}
	\| \uon\|_{L_t^2L_x^6}	\| \utw\|_{L_t^2L_x^6}
		\|\uth  \|_{L_t^\infty L_x^2}        
	\|v \|_{L_t^\infty L_x^2}.
\end{align*}
By \eqref{low str} and the embeddings $U_m^2,V_m^2\hookrightarrow L_t^\infty L_x^2$ in Lemma~\ref{Lem:property UV}, we obtain
\begin{align*}
	I&\lesssim_{N_0}
	\sum_{M\ls N_0}   \| V_M \|_{L_x^\frac32}
	\|u_1 \|_{\ul^{2,1}} 
	\|u_2 \|_{\ul^{2,1}} 
	\|\uth  \|_{L_t^\infty L_x^2}        
	\|v \|_{L_t^\infty L_x^2}.  
\end{align*}
Here, bold $\mathbf{u_3}$ means that the estimates hold for both $u_3$ and $\nabla_{\mathbb{S}}u_3$ cases.
From \eqref{ineq:potential} we estimate
$$\sum_{M\ls N_0}  \| V_M \|_{L_x^\frac32}\le\sum_{0<M \le 1} M^{1-\gamma_1}+\sum_{1< M\ls N_0} M^{1-\gamma_2}\ls C(N_0),$$ where  the assumption $\gamma_1$ be less than 1 is essential. Thus we have 
\begin{align*}
	I&\lesssim_{N_0}
	\|u_1 \|_{\ul^{2,1}} 
	\|u_2 \|_{\ul^{2,1}} 
	\|u_3  \|_{\ul^{2,1}}        
	\|v \|_{\vl^2}. 
\end{align*}

\subsubsection{ $N_0=\min(N_1,N_2)\ll \max(N_1,N_2)$ }
In this case $M$ should be comparable to $\max(N_1,N_2)$.
We divide the case according to whether $u_3$ takes the angular derivative or not.

{ \it (1) $u_3$ case:}
In this case at most one of $u_1$, $u_2$ could take the angular derivative. We denote this by bold $\uon$, $\utw$.
We have by H\"older inequality
\begin{align*}
I&\ls \sum_{M\sim \max(N_1,N_2)} \|V_M*(\uon\utw)\|_{L_t^1 \mathcal{L}_\rho^\infty L_\theta^\frac{6r}{2r-6}} 
\|\widetilde{\dot P_M} (u_3 v) \|_{L_t^\infty \mathcal{L}_\rho^1 L_\theta^\frac{6r}{4r+6}}.
\end{align*}
We compute the first norm. We assume $N_0=N_1<N_2$. We apply Lemma~\ref{radial2}
\begin{align}\label{es1}
\|V_M*(u_1u_2)\|_{L_t^1 \mathcal{L}_\rho^\infty L_\theta^\frac{6r}{2r-6}}
\lesssim \|V_M\|_{L_x^{\frac{6r}{5r-6}}} \|u_1 u_2\|_{L_t^1 \mathcal{L}_\rho^\frac{6r}{r+6} L_\theta^{2} }.
\end{align}
We estimate using Lemma~\ref{sob-s2}
\begin{align}\label{es2}
\begin{aligned}
\|u_1 u_2\|_{L_t^1 \mathcal{L}_\rho^\frac{6r}{r+6} L_\theta^{2} } 
&\lesssim  \|u_1\|_{L_t^2 \mathcal{L}_\rho^6 L_\theta^\infty}\|u_2\|_{L_t^2 \mathcal{L}_\rho^r L_\theta^2} 
\ls (\|u_1\|_{L_t^2 L_x^6} +\|\nabla_{\mathbb S}u_1\|_{L_t^2L_x^6})\|u_2\|_{L_t^2 \mathcal{L}_\rho^r L_\theta^2}  \\
&\ls_{N_0} (\|u_1\|_{U_m^2} + \|\nabla_{\mathbb S} u_1\|_{U_m^2})
N_2^\frac1r \|u_2\|_{U_m^2},
\end{aligned}
\end{align}
where in the last inequality we used Strichartz estimates \eqref{low str} and \eqref{ang str}. 
Similarly we estimate
\begin{align}\label{es3}
\begin{aligned}
\|u_1 u_2\|_{L_t^1 \mathcal{L}_\rho^\frac{6r}{r+6} L_\theta^{2} } 
&\lesssim  \|u_1\|_{L_t^2 L_x^6}\|u_2\|_{L_t^2 \mathcal{L}_\rho^r L_\theta^3} 
\ls \|u_1\|_{L_t^2 L_x^6} ( \|u_2\|_{L_t^2 \mathcal{L}_\rho^r L_\theta^2} +\|\nabla_{\mathbb S} u_2\|_{L_t^2 \mathcal{L}_\rho^r L_\theta^2}) \\
&\ls_{N_0} \|u_1\|_{U_m^2} N_2^\frac1r ( \|u_2\|_{U_m^2} + \|\nabla_{\mathbb S} u_1\|_{U_m^2}).
\end{aligned}
\end{align}
The other case $N_0=N_2<N_1$ can be bounded similarly.
Thus from \eqref{es1},\eqref{es2} and \eqref{es3} we obtain
\begin{align*}
\|V_M*(\uon\utw)\|_{L_t^1 \mathcal{L}_\rho^\infty L_\theta^\frac{6r}{2r-6}}
&\ls \|V_M\|_{L_x^{\frac{6r}{5r-6}}} \max(N_1,N_2)^\frac1r\|u_1\|_{U_m^{2,1}} \|u_2\|_{U_m^{2,1}}.
\end{align*}
Next we estimate the second norm using Lemma~\ref{radial2} 
\begin{align*}
\|\widetilde{\dot P_M} (u_3 v) \|_{L_t^\infty \mathcal{L}_\rho^1 L_\theta^\frac{6r}{4r+6}}
&\ls \| \mathcal{F}^{-1}\chi_M \|_{L_x^{1}} \|u_3\|_{L_t^\infty \mathcal L_\rho^2 L_\theta^\frac{6r}{r+6}}
\|v\|_{L_t^\infty L_x^2 } \\
&\ls ( \|u_3\|_{L_t^\infty L_x^2} +\|\nabla_{\mathbb S} u_3\|_{L_t^\infty L_x^2})\|v\|_{L_t^\infty L_x^2 },
\end{align*}
where we applied Lemma~\ref{sob-s2} with $r>\frac{10}{3}$.
In conclusion, we have 
\begin{align*}
I&\ls_{N_0} \sum_{M \sim \max(N_1,N_2)} \|V_M\|_{L_x^{\frac{6r}{5r-6}}}  \| \mathcal{F}^{-1}\chi_M \|_{L_x^{1}}
\max(N_1,N_2)^\frac1r \|u_1\|_{U_m^{2,1}} \|u_2\|_{U_m^{2,1}}
\|u_3\|_{U_m^{2,1}} \|v\|_{V_m^2} \\
&\ls \max(N_1,N_2)^\frac1r \|u_1\|_{U_m^{2,1}} \|u_2\|_{U_m^{2,1}}
\|u_3\|_{U_m^{2,1}} \|v\|_{V_m^2},
\end{align*}
where we used  $\sum_{M \sim \max(N_1,N_2)}\|V_M\|_{L_x^{\frac{6r}{5r-6}}} \| \mathcal{F}^{-1}\chi_M \|_{L_x^{1}} \ls \sum_{M\sim \max(N_1,N_2) } M^{3(\frac16+\frac1r)-\gamma_2}<C$ for $r$ we consider by \eqref{ineq:potential}.

{ \it (2) $\nabla_{\mathbb S} u_3$ case:}
In this case neither $u_1$ nor $u_2$ takes the angular derivative. We have by H\"older inequality
\begin{align}\label{22}
I&\ls \sum_{M\sim \max(N_1,N_2)} \|V_M*(u_1u_2)\|_{L_t^1 L_x^\infty } 
\|\widetilde{\dot P_M} (\nabla_{\mathbb S} u_3 v) \|_{L_t^\infty L_x^1}.
\end{align}
We consider the former. By symmetry we may assume $N_0=N_1<N_2$. We apply Lemma~\ref{sob-s2}  
\begin{align*}
\|V_M*(u_1u_2)\|_{L_t^1 L_x^\infty}
\ls  \|V_M*(u_1u_2)\|_{L_t^1 L_x^\infty L_\theta^{\frac{2r}{r-2}}}+ \|\nabla_{\mathbb S} V_M*(u_1u_2)\|_{L_t^1 L_x^\infty L_\theta^{\frac{2r}{r-2}}}.
\end{align*}
By applying Lemma~\ref{radial2} and H\"older inequality we estimate
\begin{align*}
\|V_M*(u_1u_2)\|_{L_t^1 L_x^\infty L_\theta^{\frac{2r}{r-2}}}
&\ls \|V_M\|_{L_x^{\frac{6r}{5r-6}}}\|u_1 u_2\|_{L_t^1\mathcal{L}_{\rho}^{\frac{6r}{6+r}}L_\theta^\frac32}
\ls \|V_M\|_{L_x^{\frac{6r}{5r-6}}}\| u_1\|_{L_t^2 L_x^6 } \|u_2\|_{L_t^2 \mathcal{L}_\rho^{r}L_\theta^2} \\
&\ls_{N_0} \|V_M\|_{L_x^{\frac{6r}{5r-6}}} \|u_1\|_{U_m^2} N_2^\frac1r \|u_2\|_{U_m^2},
\end{align*}
 The derivative term can be estimated by the same argument as above because by Lemma~\ref{radial1} and Leibniz rule, we have
\begin{align*}
\|\nabla_{\mathbb S} V_M*(u_1u_2)\|_{L_t^1 L_x^\infty L_\theta^{\frac{2r}{r-2}}}
\le \|V_M*(\nabla_{\mathbb S} u_1 u_2)\|_{L_t^1 L_x^\infty L_\theta^{\frac{2r}{r-2}}} + \|V_M*( u_1 \nabla_{\mathbb S} u_2)\|_{L_t^1 L_x^\infty L_\theta^{\frac{2r}{r-2}}}.
\end{align*}
Then we finally obtain
\begin{align*}
\|V_M*(u_1u_2)\|_{L_t^1 L_x^\infty}
\ls_{N_0} \|V_M\|_{L_x^{\frac{6r}{5r-6}}} \max(N_1,N_2)^\frac1r \|u_1\|_{U_m^{2,1}}\|u_2\|_{U_m^{2,1}}.
\end{align*}
For the latter in \eqref{22} we only use H\"older inequality and the embedding
\begin{align*}
\|\widetilde{\dot P_M} (\nabla_{\mathbb S} u_3 v) \|_{L_t^\infty L_x^1}
\ls \|\nabla_{\mathbb S} u_3\|_{U_m^2} \|v\|_{V_m^2}.
\end{align*}
In conclusion we get as in the previous case
\begin{align*}
I&\ls_{N_0} \sum_{M\sim \max(N_1,N_2)} M^{3(\frac16+\frac1r)-\gamma_2}\max(N_1,N_2)^\frac1r
\|u_1\|_{U_m^{2,1}}\|u_2\|_{U_m^{2,1}}\|u_3\|_{U_m^{2,1}}\|v\|_{V_m^2}\\
&\ls \max(N_1,N_2)^\frac1r
\|u_1\|_{U_m^{2,1}}\|u_2\|_{U_m^{2,1}}\|u_3\|_{U_m^{2,1}}\|v\|_{V_m^2}.
\end{align*}

\subsubsection{ $N_0<N_1,N_2$ }
We apply H\"older inequality  
\begin{align*}
I &\lesssim
\sum_{M>0} \|V_M *(\uon\utw)\|_{L_t^1 L_x^\infty} 
\|\widetilde{\dot{P}_M} (\uth v) \|_{L_t^\infty L_x^1} \\
&\lesssim
\sum_{M>0}  \|V_M\|_{L_x^{\frac{r}{r-2}}} 
\|\uon\utw\|_{L_t^1 L_x^\frac r2}  
\|\uth  \|_{L_t^\infty L_x^2}        
\|v \|_{L_t^\infty L_x^2}.        
\end{align*}
We claim that $ \|\uon\utw\|_{L_t^1 L_x^\frac r2}   \ls N_1^{\frac1r}N_2^{\frac1r}
\|u_1\|_{\ul^{2,1} }  \|u_2\|_{\ul^{2,1} }$. Indeed, we estimate applying Lemma~\ref{sob-s2} and spherical Strichartz estimate \eqref{ang str}
\begin{align*}
\|u_1u_2\|_{L_t^1 L_x^{\frac{r}{2}} }
&\ls \|u_1\|_{L_t^2 \mathcal L_\rho^r L_\theta^2 } 
\|u_2\|_{L_t^2 \mathcal L_\rho^r L_\theta^{\frac{2r}{4-r}} } \\
&\ls \|u_1\|_{L_t^2 \mathcal L_\rho^r L_\theta^2 } 
( \|u_2\|_{L_t^2 \mathcal L_\rho^r L_\theta^2 } 
+ \|\nabla_{\mathbb S} u_2\|_{L_t^2 \mathcal L_\rho^r L_\theta^2 } ) \\
&\ls N_1^{\frac1r}N_2^{\frac1r}\|u_1\|_{\ul^{2} } \big( \|u_2\|_{\ul^{2} }+ \|\nabla_{\mathbb S} u_2\|_{\ul^{2} } \big).
\end{align*}
Also, we can change the role of $u_1$ and $u_2$, which implies the claim. Thus we have
\begin{align*}
I \lesssim
\sum_{M>0}  \|V_M\|_{L_x^{\frac{r}{r-2}}} 
N_1^{\frac1r}N_2^{\frac1r}
\|u_1\|_{\ul^{2,1} }  \|u_2\|_{\ul^{2,1} } \|u_3\|_{\ul^{2,1} }
\|v \|_{L_t^\infty L_x^2}
\end{align*}
We compute the summation over $M$ using \eqref{ineq:potential}
$$ \sum_{M>0}   \|V_M\|_{L_x^{\frac{r}{r-2}}}=\sum_{0<M\le 1}M^{\frac6r-\gamma_1}+\sum_{M>1}M^{\frac6r-\gamma_2} <C, $$ 
which is finite if we choose $r$ so that $r>6/\gamma_2$.

\subsection{Case2: ${N_3\ll N}$}
In this subsection we prove
\begin{align*}
I\ls \min(N_1,N_2)^\frac1r N_3^\frac1r \|u_1\|_{\ul^{2,1} }  \|u_2\|_{\ul^{2,1} }\|u_3\|_{\ul^{2,1} }  \|v\|_{V_m^{2} }.
\end{align*}
In this case we should further divide the case whether $N_3$ is $N_0$ or not. Among them the case $N_0=\min(N_1,N_2)\sim \max(N_1,N_3)$ is already considered in section~4.1.1. 

Note that in this range we have $M\sim N\ls \max(N_1,N_2)$.

\subsubsection{ $N_0=N_3\ll N$}
Suppose $N_0=\min(N_1,N_2)\ll\max(N_1,N_2)$.
We estimate
\begin{align*}
I&\ls \sum_{M \sim N}\|V_M*(\uon\utw)\|_{L_t^2 L_x^3} 
\|\widetilde{\dot{P}_M} (\uth v) \|_{L_t^2 L_x^\frac32} \\
&\ls  \sum_{M \sim N}\|V_M\|_{L_x^\frac32} \|\uon\|_{L_t^2 L_x^6} \|\utw \|_{L_t^\infty L_x^2} \|\uth\|_{L_t^2 L_x^6} \|v \|_{L_t^\infty L_x^2} \\
&\ls_{N_0}  \sum_{M \sim N} M^{1-\gamma_2} \|u_1\|_{\ul^{2,1} }  \|u_2\|_{\ul^{2,1} }\|u_3\|_{\ul^{2,1} }  \|v\|_{V_m^{2} },
\end{align*}
which is complete since $\gamma_2>\frac32$.

Suppose $N_1,N_2>N_0$. We have
\begin{align}\label{rrr-1}
I\ls \sum_{M \sim N}\|V_M*(\uon\utw)\|_{L_t^2 \mathcal{L}_\rho^2 L_\theta^r} 
\|\widetilde{\dot{P}_M} (\uth v) \|_{L_t^2 \mathcal{L}_\rho^2 L_\theta^{\frac{r}{r-1}}}. 
\end{align}

We bound the first term. We assume $\min(N_1,N_2)=N_1$. By Lemma~\ref{radial2} we have
\begin{align}\label{22r}
\| V_M * (u_1 u_2) \|_{L_t^2 \mathcal{L}_\rho^2 L_\theta^r}
&\ls \|V_M\|_{L_x^\frac{r}{r-1}} \|u_1u_2\|_{L_t^2\mathcal{L}_\rho^{\frac{2r}{r+2}}L_\theta^\frac r2}.
\end{align}
We estimate
\begin{align}\label{22r2}
\begin{aligned}
\|u_1u_2\|_{L_t^2\mathcal{L}_\rho^{\frac{2r}{r+2}}L_\theta^\frac r2} &\ls \|u_1\|_{L_t^2 \mathcal{L}_\rho^r L_\theta^{\frac{2r}{4-r}}}
\|u_2\|_{L_t^\infty \mathcal{L}_\rho^{2}L_\theta^2} \\
&\ls \Big( \|u_1\|_{L_t^2 \mathcal{L}_\rho^r L_\theta^2}
+\|\nabla_{\mathbb S}u_1\|_{L_t^2 \mathcal{L}_\rho^r L_\theta^2} \Big)\|u_2\|_{L_t^\infty L_x^2} \\
&\ls N_1^\frac1r \Big( \|u_1\|_{U_m^2}
+\|\nabla_{\mathbb S}u_1\|_{U_m^2} \Big)\|u_2\|_{U_m^2},
\end{aligned}
\end{align}
where we used Lemma~\ref{sob-s2} since $\frac{2r}{4-r}>2$.
Or, exchanging a spherical pair for H\"older inequality  we estimate
\begin{align}\label{22r3}\begin{aligned}
\|u_1u_2\|_{L_t^2\mathcal{L}_\rho^{\frac{2r}{r+2}}L_\theta^\frac r2} &\ls \|u_1\|_{L_t^2 \mathcal{L}_\rho^r L_\theta^{2}}
\|u_2\|_{L_t^\infty \mathcal{L}_\rho^{2}L_\theta^\frac{2r}{4-r}} \\
&\ls  \|u_1\|_{L_t^2 \mathcal{L}_\rho^r L_\theta^2} \Big(\|u_2\|_{L_t^\infty L_x^2} + \|\nabla_{\mathbb S} u_2\|_{L_t^\infty L_x^2}\Big) \\
&\ls  N_1^\frac1r\|u_1\|_{U_m^2} \Big(\|u_2\|_{U_m^2} + \|\nabla_{\mathbb S} u_2\|_{U_m^2}\Big).
\end{aligned}
\end{align}
Since we can change the role of $u_1$ and $u_2$, \eqref{22r},\eqref{22r2} and \eqref{22r3} imply 
\begin{align}\label{ineq:22r}
\|V_M*(\uon\utw)\|_{L_t^2 \mathcal{L}_\rho^2 L_\theta^r}
\ls \|V_M\|_{L_x^{\frac{r}{r-1}}}
\min(N_1,N_2)^{\frac1r} \|u_1\|_{\ul^{2,1} }  \|u_2\|_{\ul^{2,1} }.
\end{align}
Next we bound the second term in \eqref{rrr-1} by applying Lemma~\ref{radial2} 
\begin{align*}
\|\widetilde{\dot{P}_M} (\uth v) \|_{L_t^2 \mathcal{L}_\rho^2 L_\theta^{\frac{r}{r-1}}}
&\ls \| \mathcal{F}^{-1}\chi_M\|_{L_x^\frac65} \|\uth v\|_{L_t^2L_x^\frac32}
\ls \| \mathcal{F}^{-1}\chi_M\|_{L_x^\frac65} \|\uth \|_{L_t^2L_x^6} \|v\|_{L_t^\infty L_x^2}\\
&\ls_{N_0} \| \mathcal{F}^{-1}\chi_M\|_{L_x^\frac65} \|\uth\|_{U_m^2}\|v\|_{V_m^2}. 
\end{align*}
In conclusion we obtain
\begin{align*}
I&\ls \sum_{M \sim N} \|V_M\|_{L_x^{\frac{r}{r-1}}}\| \mathcal{F}^{-1}\chi_M\|_{L_x^\frac65} \min(N_1,N_2)^{\frac1r} \|u_1\|_{\ul^{2,1} }  \|u_2\|_{\ul^{2,1} }
\| u_3\|_{U_m^{2,1}} \|v\|_{V_m^{2}}.
\end{align*}
which implies the desired result since we have from \eqref{ineq:potential} $$\sum_{M\sim N} \|V_M\|_{L_x^{\frac{r}{r-1}}} \|\mathcal{F}^{-1}\chi_M\|_{L_x^\frac65}  \ls \sum_{M \sim N} M^{\frac3r+\frac12-\gamma_2}<C.$$

\subsubsection{ $N_0<N_3\ll N \ \text{and}\  N_0=\min(N_1,N_2)\ll \max(N_1,N_2)  $ }
We divide the case according to whether $u_3$ takes the angular derivative or not.

{ \it (1) $u_3$ case:} We have
\begin{align*}
I\ls \sum_{M \sim N}\|V_M*(\uon\utw)\|_{L_t^2 L_x^\frac{2r}{r-2}} \|\widetilde{\dot P_M} ( u_3 v) \|_{L_t^2 L_x^{\frac{2r}{r+2}}}.
\end{align*}
We compute the first norm. By symmetry we may assume $\min(N_1,N_2)=N_1$. 
\begin{align*}
\|V_M*(\uon\utw)\|_{L_t^2 L_x^\frac{2r}{r-2}}
&\ls \|V_M\|_{L^\frac{6r}{5r-6}}\|\uon\utw\|_{L_t^2 L_x^{\frac32}}
\ls \|V_M\|_{L^\frac{6r}{5r-6}} \|\uon\|_{L_t^2L_x^6} \|\utw\|_{L_t^\infty L_x^2} \\
&\ls_{N_0} \|V_M\|_{L^\frac{6r}{5r-6}} \|\uon\|_{U_m^2} \|\utw\|_{U_m^2}.
\end{align*}
And we estimate the second term using Lemma~\ref{sob-s2}
\begin{align*}
\|\widetilde{\dot P_M} ( u_3 v) \|_{L_t^2 L_x^{\frac{2r}{r+2}}} 
&\ls \| \mathcal{F}^{-1}\chi_M \|_{L_x^{1}} \|u_3\|_{L_t^2 L_x^r} \|v\|_{L_t^\infty L_x^2}
\ls \big( \|u_3\|_{L_t^2\mathcal L_\rho^r L_\theta^2} + \| \nabla_{\mathbb S} u_3\|_{L_t^2\mathcal L_\rho^r L_\theta^2}\big)\|v\|_{L_t^\infty L_x^2}\\
&\ls N_3^\frac1r \big( \|u_3\|_{U_m^2} + \|\nabla_{\mathbb S} u_3\|_{U_m^2}\big) \|v\|_{V_m^2}.
\end{align*}
Thus we have
\begin{align*}
I&\ls \sum_{M \sim N} M^{3(\frac1r+\frac16)-\gamma_2} N_3^\frac1r \|u_1\|_{\ul^{2,1} }  \|u_2\|_{\ul^{2,1} }
\| u_3\|_{U_m^{2,1}} \|v\|_{V_m^{2}} \\
&\ls N_3^\frac1r \|u_1\|_{\ul^{2,1} }  \|u_2\|_{\ul^{2,1} }
\| u_3\|_{U_m^{2,1}} \|v\|_{V_m^{2}}.
\end{align*}

{ \it (2) $\nabla_{\mathbb S} u_3$ case:}
We have
\begin{align*}
I&\ls \sum_{M\sim \max(N_1,N_2)} \|V_M*(u_1u_2)\|_{L_t^2 \mathcal{L}_\rho^\frac{2r}{r-2}L_\theta^\infty} 
\|\widetilde{\dot P_M} (\nabla_{\mathbb S}u_3 v) \|_{L_t^2 \mathcal{L}_\rho^{\frac{2r}{r+2}}L_\theta^1}.
\end{align*}
We estimate the first norm. Applying Lemma~\ref{sob-s2} we obtain	
\begin{align*}
\|V_M*(u_1u_2)\|_{L_t^2 \mathcal{L}_\rho^\frac{2r}{r-2}L_\theta^\infty}
\ls \|V_M*(u_1u_2)\|_{L_t^2 \mathcal{L}_\rho^\frac{2r}{r-2}L_\theta^\frac{2r}{r-2}}
+ 	\| \nabla_{\mathbb S} V_M*(u_1u_2)\|_{L_t^2 \mathcal{L}_\rho^\frac{2r}{r-2}L_\theta^\frac{2r}{r-2}}
\end{align*}	
By Young's and H\"older inequality we have
\begin{align*}
\|V_M*(u_1u_2)\|_{L_t^2 \mathcal{L}_x^\frac{2r}{r-2}}
\ls \|V_M\|_{L_x^{\frac{6r}{5r-6}}} \|u_1\|_{L_t^2 L_x^6} \|u_2\|_{L_t^\infty L_x^2}
\ls_{N_0} \|V_M\|_{L_x^{\frac{6r}{5r-6}}} \|u_1\|_{U_m^2} \|u_2\|_{U_m^2}.
\end{align*}
Then by Leibniz rule we can bound the derivative term similarly and finally get 
\begin{align}\label{1}
\|V_M*(u_1u_2)\|_{L_t^2 \mathcal{L}_\rho^\frac{2r}{r-2}L_\theta^\infty}
\ls \|V_M\|_{L_x^{\frac{6r}{5r-6}}} \|u_1\|_{U_m^{2,1}}\|u_2\|_{U_m^{2,1}}.
\end{align}
We apply the Lemma~\ref{radial2} to the second term
\begin{align}\label{2}
\|\widetilde{\dot P_M} (\nabla_{\mathbb S} u_3 v) \|_{L_t^2 \mathcal{L}_\rho^{\frac{2r}{r+2}}L_\theta^1}
\ls \| \mathcal{F}^{-1}\chi_M \|_{L_x^{1}}  
\| \nabla_{\mathbb S} u_3\|_{L_t^2 \mathcal{L}_\rho^r L_\theta^2} \|v\|_{L_t^\infty L_x^2}
\ls  N_3^\frac1r 
\| \nabla_{\mathbb S} u_3 \|_{U_m^{2}} \|u_2\|_{V_m^{2}}.
\end{align} 
By \eqref{1} and \eqref{2} we obtain
\begin{align*}
I &\ls_{N_0} \sum_{M\sim N }\|V_M\|_{L_x^{\frac{6r}{5r-6}}} \|u_1\|_{\ul^{2,1} }  \|u_2\|_{\ul^{2,1}}
N_3^\frac1r 
\| u_3\|_{U_m^{2,1}} \|v\|_{V_m^{2}} \\
&\ls N_3^\frac1r \|u_1\|_{\ul^{2,1} }  \|u_2\|_{\ul^{2,1}} \| u_3\|_{U_m^{2,1}} \|v\|_{V_m^{2}}.
\end{align*}

\subsubsection{ $N_0<N_1,N_2,N_3 \ \text{and}\ N_3\ll N$ }
\begin{align*}
I &\lesssim
\sum_{M\sim N} \|V_M *(\uon\utw)\|_{L_t^2 \mathcal{L}_\rho^2 L_\theta^r} 
\|\widetilde{\dot{P}_M}(\uth v) \|_{L_t^2 \mathcal{L}_\rho^2 L_\theta^{\frac{r}{r-1}}}. 
\end{align*}
The first term is bounded as in \eqref{ineq:22r}.
For the second one we apply Lemma~\ref{radial2}
\begin{align*}
\|\widetilde{\dot{P}_M}(\uth v) \|_{L_t^2 \mathcal{L}_\rho^2 L_\theta^\frac{r}{r-1}}
&\ls \| \mathcal{F}^{-1} \chi_M \|_{L_x^{\frac{r}{r-1}}} 
\|\uth v\|_{L_t^2\mathcal{L}_\rho^{\frac{2r}{r+2}}L_\theta^1} 
\ls \| \mathcal{F}^{-1} \chi_M \|_{L_x^{\frac{r}{r-1}}} 
\|\uth\|_{L_t^2 \mathcal{L}_\rho^{r}L_\theta^2}
\|v\|_{L_t^\infty \mathcal{L}_\rho^{2}L_\theta^2} \\
&\ls \| \mathcal{F}^{-1} \chi_M \|_{L_x^{\frac{r}{r-1}}} 
N_3^\frac1r \|\uth\|_{U_m^2}\|v\|_{V_m^2} .
\end{align*}
Then the claim follows since $\sum_{M\sim N} \|V_M\|_{L_x^{\frac{r}{r-1}}}\| \mathcal{F}^{-1} \chi_M \|_{L_x^{\frac{r}{r-1}}} 
\ls \sum_{M\sim N} M^{\frac6r-\gamma_2} < C $ by \eqref{ineq:potential}.



\section*{Acknowledgements}
The author would like to thank Prof. Yonggeun Cho for his encouragement and advice on the paper. And the author is grateful to the referee for careful reading of the paper and valuable comments. 
This work was supported by NRF (NRF-2015R1D1A1A09057795).



\medskip


\begin{thebibliography}{00}

\bibitem{BH2017} I. Bejenaru and S. Herr
{\it On global well-posedness and scattering for the massive
Dirac-Klein-Gordon system}, J. Eur. Math. Soc. \textbf{19} (2017), no. 8, 2445--2467.

\bibitem{cg 76} Chadam, J. M. and  Glassey, R. T. ,
{\it On the Maxwell-Dirac equations with zero magnetic field and their solution in two space dimensions}, J. Math. Anal. Appl. \textbf{53}
(1976), no. 3, 495--507.

\bibitem{chonak} Y. Cho and K. Nakanishi, {\it On the global existence of semirelativistic Hartree equations}, RIMS Kokyuroku Bessatsu, \textbf{B22} (2010), 145-166.

\bibitem{choz0}  Y. Cho and T. Ozawa,{\it On the semi-relativisitc Hartree type equation}, SIAM J. Math. Anal.,
\textbf{38} (2006), no. 4, 1060--1074.

\bibitem{co07} Y. Cho and T. Ozawa,
{\it Global solutions of semirelativistic Hartree type equations}, J. Korean Math. Soc., \textbf{44} (2007), no.5, 1065--1078.

\bibitem{coss 09} Y. Cho, T. Ozawa, H. Sasaki and Y. Shim {\it Remarks on the semirelativistic Hartree equations}, Discrete Contin. Dyn. Syst.
\textbf{ 23} (2009), no. 4, 1277--1294.

\bibitem{cox}  Y. Cho, T. Ozawa and S. Xia, {\it Remarks on some dispersive estimates}, Commun. Pure Appl. Anal., {\bf 10} (2011), no. 4,
1121-1128.


\bibitem{crt-n} T. Coulhon, E. Russ and V. Tardivel-Nachef, {\it Sobolev Algebras on Lie Groups and Riemannian Manifolds}, American Journal of Mathematics \textbf{123}, no. 2 (2001), 283--342.

\bibitem{dfs07}
P. D'Ancona, D. Foschi and S. Selberg, 
{\it Null structure and almost optimal local regularity for the
	Dirac-Klein-Gordon system}, J. Eur. Math. Soc. (JEMS)
\textbf{9} (2007), no.4, 877--899.

\bibitem{E07} A. Elgart,  {\it Mean field dynamics of boson stars},
Comm. Pure Appl. Math. \textbf{60}, No. 4 (2007) 500--545.

\bibitem{l03}  J. Fr\"ohlich and E. Lenzmann,  {\it Mean-field limit of quantum {B}ose gases and nonlinear
	{H}artree equation}, S\'eminaire: \'Equations aux {D}\'eriv\'ees {P}artielles. 2003--2004  \textbf{S\'emin. \'Equ. D\'eriv. Partielles} (2004), Exp. no. XIX, 26.


\bibitem{ghn} Z. Guo, Z. Hani and K. Nakanishi, {\it Scattering for the 3D Gross-Pitaevskii equation}, Comm. Math. Phys. \textbf{259} (2018), no.1, 265-295.

\bibitem{HL2014} S. Herr and E. Lenzmann, 
{\it The Boson star equation with initial data of low regularity}, Nonlinear Anal. \textbf{97} (2014), 125--137.

\bibitem{hhk} M. Hadac, S. Herr and H. Koch, {\it Well-posedness and scattering for the KP-II equation in a critical space}, Ann. Inst. H. Poincare Anal. Non Lineaire  \textbf{26}  (2009), no. 3, 917--941.


\bibitem{hete} S. Herr and T. Tesfahun, {\it Small data scattering for semi-relativistic equations with Hartree type nonlinearity}, J. Differential Equations  \textbf{259}, (2015), no.10, 5510--5532.


\bibitem{ktv} H. Koch, D. Tataru and M. Visan, {\it Dispersive Equations and Nonlinear Waves}, Oberwolfach Seminars \textbf{45} (2014).


\bibitem{l07} E. Lenzmann, {\it Well-posedness for semi-relativistic Hartree equations of critical type}, Math. Phys. Anal. Geom. \textbf{10} (2007), no. 1, 43--64.

\bibitem{EH87} E, Lieb and H. Yau, {\it The Chandrasekhar theory of stellar collapse as the limit of quantum mechanics}, Comm. Math. Phys. \textbf{112} (1987), no.1, 147--174.

\bibitem{F2014} F. Pusateri, {\it Modified scattering for the boson star equation}, Comm. Math. Phys. \textbf{332} (2014), no.3, 1203--1234.


\end{thebibliography}
\end{document}